\documentclass[reqno]{amsart}
\usepackage{graphicx}

\theoremstyle{plain}
\newtheorem{theorem}{Theorem}
\newtheorem{lemma}{Lemma}

\theoremstyle{definition}
\newtheorem{definition}{Definition}

\theoremstyle{remark}
\newtheorem{remark}{Remark}

\title[Singular shock waves in interactions]
{Singular shock waves in interactions
\footnote{\tiny  The work is supported by Serbian Ministry of 
Science and Ecology} }
\author{Marko Nedeljkov} 
\address{Department of Mathematics and Informatics, University of Novi Sad,
Trg D. Obradovi\'{c}a 4, 21000 Novi Sad, Yugoslavia} 
\email{markonne@uns.ns.ac.yu}
\date{}

\begin{document}

\begin{abstract}
In a number of papers it was shown that there are one-dimensional
systems such that they contain solutions with, so called, 
overcompressive singular shock waves
besides the usual elementary waves (shock and rarefaction ones
as well as contact discontinuities). 

One can see their definition for a general 2 $\times$ 2 system
with fluxes linear in one of dependent variables in \cite{Ned1}.
This paper is devoted to examining 
their interactions with themselves and elementary waves.
After a discussion of systems given in a general form, a complete
analysis will be given for the ion-acoustic system given in \cite{KeyKr}.
\medskip


\end{abstract}

\maketitle
\section{Introduction}

Consider the system
\begin{equation}\label{gdss1}
\begin{split}
& (f_{2}(u))_{t}+(f_{3}(u)v+f_{4}(u))_{x}=0 \\
& (g_{1}(u)v+g_{2}(u))_{t}+(g_{3}(u)v+g_{4}(u))_{x}=0.
\end{split}
\end{equation}
where $f_{i},g_{j}$, $i=2,...,4$, $j=1,...,4$ are polynomials with the 
maximal degree $m$, $(u,v)=(u(x,t),v(x,t))$ are unknown functions with 
a physical range $\Omega$, $(x,t)\in {\mathbb R}\times {\mathbb R}_{+}$. 
We shall fix the following notation for the rest of the paper:
$$f_{i}(y)=\sum_{k=0}^{m}a_{i,k}y^{k}, \; 
g_{j}(y)=\sum_{k=0}^{m}b_{j,k}y^{k}, \; i=2,3,4,\; j=1,2,3,4. $$

There are cases when there is no classical solution to Riemann problem
for the above system. Sometimes, there is a solution in the form
of delta or singular shock wave. 
In \cite{Ned1} one can see when 
a system in evolution form (i.e.\ 
when $f_{2}=u$, $g_{1}=1$ and $g_{2}=0$) permits a solution in the 
shape of singular shock wave. With the same type of reasoning
and a more effort, one can give the answer to the 
same question in the case system (\ref{gdss1}).

The aim of this paper is to investigate what happens during and after
an interaction of a singular shock wave with another wave. 
After a general statement about 
new initial data taken at interaction point
(of course, true for delta shock waves, too) in Section 3, 
we shall present a detailed
investigation in the case of the system (so called ion-acoustic system)
\begin{equation}\label{dss1}
\begin{split}
u_{t}+(u^{2}-v)_{x} & = 0 \\
v_{t}+(u^{3}/3-u)_{x} & = 0 
\end{split}
\end{equation}
given in \cite{KeyKr}.

Definitions and concepts used here are from \cite{Ned1}, based on 
the use of Colombeau generalized functions defined in \cite{ObeWa}.
They will be briefly described in the Section 2.
If one is not familiar with these concepts, he/she 
can assume that a solution to the above system is given 
by nets of smooth functions with equality substituted by
a distributional limit.
The reason why the generalized functions are used is 
to give opportunity for extending the procedure in this paper for arbitrary 
initial data when a system posses singular or delta shock wave 
as a solution.

Few interesting facts observed during the investigations of system
(\ref{dss1}) are arising a question about possibilities in a general
case. Observed facts are:
\begin{enumerate}
\item The singular shock wave solution to a Riemann problem for
(\ref{dss1}) always has an increasing strength
of the rate ${\mathcal O}(t)$, $t \to \infty$. (The strength of the 
shock is a function which multiplies the delta function contained in 
a solution, $s(t)$ in (\ref{prom})).
After the interaction, the resulting singular shock wave is 
supported by a curve, not necessary straight
line as before, and its strength can be 
an increasing, but also a constant or a decreasing function with 
the respect to the time variable.
\item When the resulting singular shock wave has a decreasing strength
(this can occurs during an interaction of a admissible 
singular shock wave with a rarefaction wave),
after some time it can decompose into two shock waves. This is a 
quite new phenomenon.
\end{enumerate}

The structure of this paper can be described in the following way.

In the second section we will introduce necessary notation and give
basic notions based on the papers \cite{ObeWa} and \cite{Ned1}.

In the third section, one can find a way how to continue a solution
to the general case of system (\ref{gdss1}) 
after an interaction point (Theorem \ref{glavna}). 
The basic assumption is that a left-hand side of the first, and
the right-hand side of the second wave can be connected by a new
singular shock wave. The conditions for such a possibility are 
formulated trough a notion of {\it second delta singular locus}, 
see Definition \ref{2d1}. Explicit calculations for a geometric description
of the locus are possible to perform for
system (\ref{gdss1}), but we shall omit it, to preserve readers
attention on the further topics.

The results given in these sections are used in the next one
devoted to special case (\ref{dss1}).

The first part of the fourth section is devoted to 
description of a situation which can occur 
after a singular shock and a shock wave interact. In the same
way one can do the same for two singular shock waves, 
as one can at the end of this section.

The final, 5th section, contains the most interesting and important 
results about singular shock and rarefaction wave interaction. 
In that case the decoupling of a singular shock into a pair of shock waves,
already mentioned before, can occur. The analysis is done when a singular
shock wave is on the left-hand side of a rarefaction wave. But
one can easily see that these results can be obtained using the same
procedure when a singular shock is on the other side
of a rarefaction wave.

\section{Notation}

We shall briefly repeat some definitions of 
Colombeau algebra given in \cite{ObeWa} and \cite{Ned1}. Denote 
${\mathbb R}_{+}^{2}:={\mathbb R}\times (0,\infty)$, 
$\overline{{\mathbb R}_{+}^{2}}:={\mathbb R}\times [0,\infty)$ and let
$C_{b}^{\infty}(\Omega)$ be the algebra of smooth functions on 
$\Omega$ bounded together with all their derivatives. Let 
$C_{\overline{b}}^{\infty}({\mathbb R}_{+}^{2})$ be a set of 
all functions $u\in C^{\infty}({\mathbb R}_{+}^{2})$ satisfying 
$u|_{{\mathbb R}\times (0,T)} \in C_{b}^{\infty}({\mathbb R}\times (0,T))$
for every $T>0$. Let us remark that every element of 
$C_{b}^{\infty}({\mathbb R}_{+}^{2})$ has a smooth extension up to 
the line $\{ t=0\}$, i.e.\ $C_{b}^{\infty}({\mathbb R}_{+}^{2})=
C_{b}^{\infty}(\overline{{\mathbb R}_{+}^{2}})$. This is also true for
$C_{\overline{b}}^{\infty}({\mathbb R}_{+}^{2})$.

\begin{definition}\label{emn}
${\mathcal E}_{M,g}({\mathbb R}_{+}^{2})$ is the set of all maps 
$G:(0,1)\times {\mathbb R}_{+}^{2} \rightarrow {\mathbb R}$, 
$(\varepsilon,x,t) \mapsto G_{\varepsilon}(x,t)$, where  
for every $\varepsilon \in (0,1)$,
$G_{\varepsilon}\in  C_{\overline{b}}^{\infty}({\mathbb R}_{+}^{2})$
satisfies:

\noindent
For every $(\alpha,\beta)
\in {\mathbb N}_{0}^{2}$ and $T>0$, there exists $N\in {\mathbb N}$ such that
$$\sup_{(x,t)\in {\mathbb R}\times (0,T)}
|\partial_{x}^{\alpha}\partial_{t}^{\beta} G_{\varepsilon}(x,t)|
={\mathcal O}(\varepsilon^{-N}), \text{ as } \varepsilon \rightarrow 0.$$

${\mathcal E}_{M,g}({\mathbb R}_{+}^{2})$ is an multiplicative differential
algebra, i.e.\
a ring of functions with the usual operations of addition and multiplication,
and differentiation which satisfies Leibniz rule.

${\mathcal N}_{g}({\mathbb R}_{+}^{2})$ is the set of all
$G\in {\mathcal E}_{M,g}({\mathbb R}_{+}^{2})$,
satisfying:

\noindent
For every $(\alpha,\beta)
\in {\mathbb N}_{0}^{2}$, $a\in {\mathbb R}$ and $T>0$
$$\sup_{(x,t)\in {\mathbb R}\times (0,T)}
|\partial_{x}^{\alpha}\partial_{t}^{\beta} G_{\varepsilon}(x,t)|
={\mathcal O}(\varepsilon^{a}), \text{ as } \varepsilon \rightarrow 0.$$
$\Box$
\end{definition}

Clearly, ${\mathcal N}_{g}({\mathbb R}_{+}^{2})$ is 
an ideal of the multiplicative
differential algebra ${\mathcal E}_{M,g}({\mathbb R}_{+}^{2})$, i.e.\
if $G_{\varepsilon}\in {\mathcal N}_{g}({\mathbb R}_{+}^{2})$ and 
$H_{\varepsilon}\in {\mathcal E}_{M,g}({\mathbb R}_{+}^{2})$, 
then $G_{\varepsilon}H_{\varepsilon}\in {\mathcal N}_{g}({\mathbb R}_{+}^{2})$. 

\begin{definition}\label{g}
The  multiplicative differential algebra 
${\mathcal G}_{g}({\mathbb R}_{+}^{2})$ 
of generalized functions is defined by 
${\mathcal G}_{g}({\mathbb R}_{+}^{2})=
{\mathcal E}_{M,g}({\mathbb R}_{+}^{2})/{\mathcal N}_{g}({\mathbb R}_{+}^{2})$.
All operations in ${\mathcal G}_{g}({\mathbb R}_{+}^{2})$ are 
defined by the corresponding ones in ${\mathcal E}_{M,g}({\mathbb R}_{+}^{2})$.
$\Box$
\end{definition}

If $C_{b}^{\infty}({\mathbb R})$ is used instead of
$C_{b}^{\infty}({\mathbb R}_{+}^{2})$ (i.e.\
drop the dependence on the $t$ variable), then  
one obtains ${\mathcal E}_{M,g}({\mathbb R})$, ${\mathcal N}_{g}({\mathbb R})$,
and consequently, the space of generalized functions on a real line, 
${\mathcal G}_{g}({\mathbb R})$.

In the sequel, $G$ denotes an element (equivalence class) 
in ${\mathcal G}_{g}(\Omega)$ 
defined by its representative
$G_{\varepsilon}\in {\mathcal E}_{M,g}(\Omega)$.

Since $C_{\overline{b}}^{\infty}({\mathbb R}_{+}^{2})=
C_{\overline{b}}^{\infty}(\overline{{\mathbb R}_{+}^{2}})$, 
one can define a restriction of a generalized function to $\{ t=0\}$ 
in the following way.

For given $G\in {\mathcal G}_{g}({\mathbb R}_{+}^{2})$, its restriction
$G|_{t=0}\in {\mathcal G}_{g}({\mathbb R})$ is the class determined 
by a function
$G_{\varepsilon}(x,0)\in {\mathcal E}_{M,g}({\mathbb R})$. 
In the same way as above, 
$G(x-ct)\in {\mathcal G}_{g}({\mathbb R})$ 
is defined by $G_{\varepsilon}(x-ct)\in {\mathcal E}_{M,g}({\mathbb R})$. 

If $G \in {\mathcal G}_{g}$ and $f\in C^{\infty}({\mathbb R})$ 
is polynomially bounded
together with all its derivatives, then one can 
easily show that the composition $f(G)$, 
defined by a representative $f(G_{\varepsilon})$,
$G\in {\mathcal G}_{g}$ makes sense. It means that $f(G_{\varepsilon})\in 
{\mathcal E}_{M,g}$ if $G_{\varepsilon}\in {\mathcal E}_{M,g}$, and 
$f(G_{\varepsilon})-f(H_{\varepsilon}) \in {\mathcal N}_{g}$ if 
$G_{\varepsilon}-H_{\varepsilon}\in {\mathcal N}_{g}$.

The equality in the space of the generalized functions ${\mathcal G}_{g}$ is
to strong for our purpose, so we need to define
a weaker relation called association. 
\begin{definition}\label{ass}
A generalized function $G\in {\mathcal G}_{g}(\Omega)$ is 
said to be {\em associated with} $u\in {\mathcal D}'(\Omega)$, $G \approx u$,
if for some (and hence every) representative
$G_{\varepsilon}$ of $G$, $G_{\varepsilon} \rightarrow u$ in 
${\mathcal D}'(\Omega)$ as $\varepsilon \rightarrow 0$.
Two generalized functions $G$ and $H$ are said to be
associated, $G\approx H$, if $G-H \approx 0$. The rate of convergence
in ${\mathcal D}'$ with respect to $\varepsilon$
is called the {\em rate of association}.
$\Box$
\end{definition}

A generalized function $G$ is said to be {\em of a bounded type} if
$$\sup_{(x,t)\in {\mathbb R}\times (0,T)} |G_{\varepsilon}(x,t)|
={\mathcal O}(1) \text{ as } \varepsilon \rightarrow 0,$$
for every $T>0$.

Let $u \in {\mathcal D}_{L^{\infty}}'({\mathbb R})$. Let ${\mathcal A}_{0}$ 
be the  set of all functions $\phi \in C_{0}^{\infty}({\mathbb R})$ satisfying
$\phi(x)\geq 0$, $x\in {\mathbb R}$,
$\int \phi(x) dx=1$ and $\operatorname{supp}\phi \subset [-1,1]$, i.e.\
\begin{equation*}
{\mathcal A}_{0}=\{\phi\in C_{0}^{\infty}:\; 
(\forall x\in {\mathbb R}) \phi(x)\geq 0,\;
\int\phi(x)dx=1,\; \mathop{\rm supp}\phi\subset [-1,1]\}.
\end{equation*}

Let $\phi_{\varepsilon}(x)=\varepsilon^{-1}\phi(x/\varepsilon)$,
$x\in {\mathbb R}$. Then 
$$ \iota_{\phi}: u \mapsto u\ast \phi_{\varepsilon}/{\mathcal N}_{g},$$
where  $u\ast \phi_{\varepsilon}/{\mathcal N}_{g}$ denotes the 
equivalence class with respect to the ideal ${\mathcal N}_{g}$,
defines a mapping of ${\mathcal D}_{L^{\infty}}'({\mathbb R})$ into 
${\mathcal G}_{g}({\mathbb R})$, where $\ast$ denotes the usual 
convolution in ${\mathcal D}'$. It is clear that $\iota_{\phi}$ commutes
with the derivation, i.e.\
$$\partial_{x}\iota_{\phi}(u)=\iota_{\phi}(\partial_{x}u).$$

\begin{definition}\label{s-d}
\begin{itemize}
\item[(a)] $G \in {\mathcal G}_{g}({\mathbb R})$ is said to be 
{\em a generalized step function} with value $(y_{0},y_{1})$ 
if it is of bounded type and 
$$ G_{\varepsilon}(y)=
\begin{cases}
y_{0}, \;  & y< -\varepsilon \\ y_{1}, \;  & y> \varepsilon
\end{cases}
$$
Denote $[G]:=y_{1}-y_{0}$.
\item [(b)] $D \in {\mathcal G}_{g}({\mathbb R})$ 
is said to be {\em generalized split 
delta function} ({\em S$\delta$-function}, for short) with value 
$(\alpha_{0},\alpha_{1})$ if 
$D=\alpha_{0}D^{-}+\alpha_{1}D^{+}$, where $\alpha_{0}+\alpha_{1}=1$
and 
\begin{equation}\label{DG}
DG\approx (y_{0}\alpha_{0}+y_{1}\alpha_{1})\delta,
\end{equation} 
for every generalized step function $G$ with value $(y_{0},y_{1})$.
\item[(c)]Let $m$ be an odd positive integer.
A generalized function $d \in {\mathcal G}_{g}({\mathbb R})$ is said to be
{\em $m'$-singular delta function} ({\em $m'$SD-function}, for short) 
with value $(\beta_{0},\beta_{1})$ if 
$d=\beta_{0}d^{-}+\beta_{1}d^{+}$, $\beta_{0}^{m-1}+\beta_{1}^{m-1}=1$, 
$d^{\pm}\in {\mathcal G}_{g}({\mathbb R})$, 
$(d^{\pm})^{i} \approx 0$, $i\in \{1,\dots,m-2,m\}$, 
$(d^{\pm})^{m-1}\approx \delta$, and 
\begin{equation}\label{m1dG}
d^{m-1}G\approx (y_{0}\beta_{0}^{m-1}+y_{1}\beta_{1}^{m-1})\delta,
\end{equation} 
for every generalized step function $G$ with value $(y_{0},y_{1})$.
\item[(d)]Let $m$ be an odd positive integer.
A generalized function $d \in {\mathcal G}_{g}({\mathbb R})$ is said to be
{\em $m$-singular delta function} ({\em $m$SD-function}, for short) 
with value $(\beta_{0},\beta_{1})$ if 
$d=\beta_{0}d^{-}+\beta_{1}d^{+}$, $\beta_{0}^{m}+\beta_{1}^{m}=1$, 
$d^{\pm}\in {\mathcal G}_{g}({\mathbb R})$, 
$(d^{\pm})^{i} \approx 0$, $i\in \{1,\dots,m-1\}$, 
$(d^{\pm})^{m}\approx \delta$, and 
\begin{equation}\label{mdG}
d^{m}G\approx (y_{0}\beta_{0}^{m}+y_{1}\beta_{1}^{m})\delta, 
\end{equation}
for every generalized step function $G$ with value $(y_{0},y_{1})$.
\end{itemize}
$\Box$
\end{definition}

In this paper we shall assume the compatibility condition
$Dd\approx 0$, where $D$ is S$\delta$- and $d$ is $m$SD- or $m'$Sd-function.

Suppose that the initial data are given by 
\begin{equation} \label{id1}
u|_{t=T}=\begin{cases} u_{0}, & x<X \\ u_{1}, & x>X \end{cases}
\; v|_{t=T}=\begin{cases} v_{0}, & x<X \\ v_{1}, & x>X. \end{cases} 
\end{equation}

\begin{definition} \label{singsh}
{\em Singular shock wave} (DSSW for short) is an associated 
solution to (\ref{dss1}) with the initial data (\ref{id1}) of the form 
\begin{equation} \label{prom} 
\begin{split}
& u((x-X),(t-T))=G((x-X)-c(t-T)) \\
& +\tilde{s}(t)(\alpha_{0}d^{-}((x-X)-c(t-T))
+\alpha_{1}d^{+}((x-X)-c(t-T))) \\
& v((x-X),(t-T))=H((x-X)-c(t-T))\\
& +s(t)(\beta_{0}D^{-}((x-X)-c(t-T))
+\beta_{1}D^{+}((x-X)-c(t-T))) \\
& + \tilde{\tilde{s}}(t)(\gamma_{0}d^{-}((x-X)-c(t-T))
+\gamma_{1}d^{+}((x-X)-c(t-T)))
\end{split}
\end{equation} 
where 
\begin{itemize}  
\item[(i)] $c\in {\mathbb R}$ is the speed of the wave, 
\item[(ii)] $s(t)$, $\tilde{s}(t)$ and $\tilde{\tilde{s}}$
are smooth functions for $t\geq 0$, and equal zero at $t=T$.
\item[(iii)] $G$ and $H$ are generalized step functions 
with values $(u_{0},u_{1})$
and $(v_{0},v_{1})$ respectively, 
\item[(iv)] $d_{1}=\alpha_{0}d^{-}+\beta_{1}d^{+}$ and 
$d_{2}=\gamma_{0}d^{-}+\gamma_{1}d^{+}$ are 
$m$SD- or $m'$SD-functions,
\item[(v)] $D=\alpha_{0}D^{-}+\alpha_{1}D^{+}$ is an S$\delta$-function 
compatible with $d$. 
\end{itemize}

The {\it singular part} of the wave is 
$$
\left[ \begin{matrix}\tilde{s}(t)(\alpha_{0}d^{-}+\alpha_{1}d^{+}) \\
s(t)(\beta_{0}D^{-}+\beta_{1}D^{+})+
\tilde{\tilde{s}}(t)(\gamma_{0}d^{-}+\gamma_{1}d^{+}) \end{matrix}\right].
$$
The wave is {\em overcompressive} 
if its speed is less or equal to 
the left- and greater or equal to the right-hand side 
characteristics i.e.\
\begin{equation*}
\lambda_{2}(u_{0},v_{0})>\lambda_{1}(u_{0},v_{0})\geq c
\geq \lambda_{2}(u_{1},v_{1})>\lambda_{1}(u_{1},v_{1}).
\end{equation*}
$\Box$
\end{definition}

\begin{remark}\label{konstrukcija}
(a) In \cite{Ned1} one can find special choice for S$\delta$- and 
and $d$ is $m$SD- or $m'$Sd-functions. For example
$D^{\pm} \in {\mathcal G}_{g}({\mathbb R})$ 
are given by the representatives
$$D_{\varepsilon}^{\pm}(y)
:={1 \over \varepsilon} \phi\Big({y-(\pm 2\varepsilon) \over \varepsilon}\Big),
\; \phi \in {\mathcal A}_{0}.$$
$m$SD- and $m'$SD-functions can be chosen in the same manner.

\noindent
(b) Compatibility condition for an S$\delta$-function $D$ and 
an $m$SD- or $m'$SD-function $d$ is automatically fulfilled if 
\begin{equation*}
\mathop{\rm supp} d_{\varepsilon}^{+} \cap  
\mathop{\rm supp} D_{\varepsilon}^{+} =
\mathop{\rm supp} d_{\varepsilon}^{-} \cap  
\mathop{\rm supp} D_{\varepsilon}^{-} = \emptyset
\end{equation*}

\noindent
(c) Idea behind the above definition of products (\ref{DG}), (\ref{m1dG})
and (\ref{mdG}) is the following. Starting point is that we know nothing
about infinitesimal values of the initial data (carried on by step functions
$G$ and $H$ above) around zero, but only that any such unmeasurable
influence stops at the points $\pm\varepsilon$. The above mentioned
definitions are made in order to get uniqueness of all products where 
step functions, S$\delta$-, $m$SD- and $m'$SD-functions appear. 
With an additional information for $G_{\varepsilon}$ and $H_{\varepsilon}$
around zero, one can choose $D$ and $d$ much more freely. For example,
in $G_{\varepsilon}$ and $H_{\varepsilon}$ are monotone functions
(which is quite natural assumption), relation (\ref{DG}) can be substituted by
$$
DG\approx \gamma \delta, \; \gamma \text{ can be any real between }
\min\{y_{0},y_{1}\} \text{ and } \max\{y_{0},y_{1}\}.
$$
The possibilities in Colombeau algebra are even wider for specific systems
instead of general case (\ref{gdss1}). One can look in \cite{Col92} for
a good review of such possibilities.
We dealing with a system in a general form and it is the reason
for using the above definition.

\noindent
(d) Due to absence of known additional facts for 
the general case (\ref{gdss1}) (hyperbolicity, additional conservation
laws,...), one can use the overcompressibility as an admissibility
condition.
\end{remark}
\medskip

\begin{definition} \label{d1}
The set of all points $(u_{1},v_{1})\in \Omega$ such that there exists an
singular shock wave solution (called {\it corresponding DSSW})
to Cauchy problem 
(\ref{gdss1},\ref{id1}) is called {\em delta singular locus}. 
We shall write $(u_{1},v_{1})\in \mathop{\rm DSL}(u_{0},v_{0})$.
If the corresponding DSSW is overcompressive,
then it is called {\em overcompressive delta singular locus}.
We shall write $(u_{1},v_{1})\in \mathop{\rm DSL}^{\ast}(u_{0},v_{0})$.
$\Box$ 
\end{definition}

In the sequel, the term ``solution'' will denote generalized function
which solves a system in the association sense.

\section{The new initial data}

Suppose that system (\ref{gdss1}) posses a DSSW solution 
for some initial data. Assume one of the following.
\begin{itemize}
\item[(i)] If an $m$SD-function is contained in the above DSSW, then assume
\begin{equation}\label{deg}
\mathop{\rm deg}(g_{1})<m-1,\;
\mathop{\rm deg}(g_{2}) < m,\; 
\mathop{\rm deg}(f_{2}) < m.
\end{equation}
\item[(ii)] If an $m'$SD-function is contained 
in the above DSSW, then assume
\begin{equation}\label{deg'}
\mathop{\rm deg}(g_{1}) < m-2, \;  
\mathop{\rm deg}(g_{2}) < m-1,\;
\mathop{\rm deg}(f_{2}) < m-1.
\end{equation}
\end{itemize}

Take the new initial data 
\begin{equation} \label{eq6}
u|_{t=T}=\begin{cases} u_{0}, & x<X \\ u_{1}, & x>X \end{cases}, 
\; v|_{t=T}=\begin{cases} v_{0}, & x<X \\ v_{1}, & x>X \end{cases} 
+ \zeta \delta_{(X,T)},
\end{equation}
for system (\ref{gdss1}), where $\zeta$ is a non-zero real.

\begin{definition} \label{2d1}
The set of all points $(u_{1},v_{1})\in \Omega$ such that there exists an
DSSW solution (called corresponding DSSW) to Cauchy problem 
(\ref{gdss1},\ref{eq6}) for some $\zeta$ is called {\em second
delta singular locus} of initial strength $\zeta$ 
for $(u_{0},v_{0})$. We shall write 
$(u_{1},v_{1})\in \mathop{\rm SDSL}_{\zeta}(u_{0},v_{0})$
If the the corresponding DSSW is overcompressive,
then it is called {\em overcompressive second delta singular locus},
and write 
$(u_{1},v_{1})\in \mathop{\rm SDSL}_{\zeta}^{\ast}(u_{0},v_{0})$.
$\Box$ 
\end{definition}

Before the main theorem, let us give a useful lemma.

\begin{lemma} \label{podskup}
Suppose that $(u_{1},v_{1})\in \mathop{\rm DSL}(u_{0},v_{0})$.
Then $(u_{1},v_{1})\in \mathop{\rm SDSL}_{\zeta}(u_{0},v_{0})$,
if $\zeta>0$. 

If the corresponding DSSW contains $m$SD-function, and $m$
is an odd number, then the statement holds true for every real $\zeta$.

Additionally, $\beta_{i}$, $i=1,2$, from Definition \ref{singsh} for the
corresponding DSSW do not depend on $\zeta$.
\end{lemma} 
\begin{proof}
We shall give the proof for a DSSW containing $m$SD-function (\ref{prom}).
The other case can be proved in the same way.

Inserting functions $u$ and $v$ from (\ref{prom}) into system (\ref{gdss1})
with initial data (\ref{eq6}) and taking account relations
(\ref{deg}) or (\ref{deg'}), one gets
\begin{equation*}
\begin{split}
f_{2}(u) \approx & f_{2}(G) \\
g_{1}(u) \approx & g_{1}(G) \\
g_{2}(u) \approx & g_{2}(G) \\
f_{3}(u) \approx & f_{3}(G) 
+ \tilde{s}(t)^{m-1}(u_{1}\alpha_{0}^{m-1}d^{-}+u_{0}\alpha_{1}^{m-1}d^{+})
m a_{3,m-1} \\
& + \tilde{s}(t)^{m}(\alpha_{0}^{m}d^{-}+\alpha_{1}^{m}d^{+}) a_{3,m} 
\tilde{s}(t)^{m}a_{3,m}\delta \\
f_{4}(u)\approx & f_{4}(G)
+ \tilde{s}(t)^{m}(\alpha_{0}^{m-1}d^{+}+\alpha_{1}^{m-1}d^{+}) a_{4,m} 
\tilde{s}(t)^{m}a_{4,m}\delta \\
g_{3}(u) \approx & g_{3}(G) 
+ s(t)^{m-1}(u_{1}\beta_{0}^{m-1}d^{-}+u_{0}\beta_{1}^{m-1}d^{+})
m b_{3,m-1} \\
& + \tilde{s}(t)^{m}(\beta_{0}^{m}d^{-}+\beta_{1}^{m}d^{+}) b_{3,m} 
\tilde{s}(t)^{m}b_{3,m}\delta \\
g_{4}(u)\approx & f_{4}(G)
+ \tilde{s}(t)^{m}(\beta_{0}^{m-1}d^{+}+\beta_{1}^{m-1}d^{+}) b_{4,m}
\tilde{s}(t)^{m}b_{4,m}\delta 
\end{split}
\end{equation*}
There are two possible cases. Either $\tilde{\tilde{s}}\not \equiv 0$ and
$a_{3,m}=b_{3,m}=0$ (i.e.\ 
$\mathop{\rm deg}(f_{3})\leq m-1$ and $\mathop{\rm deg}(g_{3})\leq m-1$),
or $\tilde{\tilde{s}}\equiv 0$. In both the cases, the procedure which 
follows is the same, so take $\tilde{\tilde{s}}\not \equiv 0$ 
for definiteness. From the first equation of (\ref{gdss1}) one gets
\begin{equation*}
\begin{split}
& (f_{2}(u))_{t}+(f_{3}(u)v+f_{4}(u))_{x}  \\
\approx & -c([f_{2}(G)]+[f_{3}(G)H+f_{4}(G)])\delta \\
& + \tilde{s}(t)^{m-1}\tilde{\tilde{s}}(t)
(u_{1}\alpha_{0}^{m-1}\gamma_{0}+u_{0}\alpha_{1}^{m-1}\gamma_{1})
ma_{3,m-1}\delta' \\
& + (f_{3}(u_{0})\beta_{0}+f_{3}(u_{1})\beta_{1})\delta'
+ \tilde{s}(t)^{m}\delta' \approx 0. 
\end{split}
\end{equation*}
One immediately gets the speed of DSSW,
$$ c={[f_{3}(G)H+f_{4}(G)] \over [f_{2}(G)]}, $$
and the relations
$$ \kappa_{1}s(t)=\tilde{s}(t)^{m-1}\tilde{\tilde{s}}(t)
\text{ and } 
\kappa_{2}s(t)=\tilde{s}(t)^{m},$$
for some reals $\kappa_{1}$ and $\kappa_{2}$. Finally, one gets
\begin{equation}\label{esp1} 
\kappa_{1}(u_{1}\alpha_{0}^{m-1}\gamma_{0}+u_{0}\alpha_{1}^{m-1}\gamma_{1})
m a_{3,m-1}+f_{3}(u_{0})\beta_{0}+f_{3}(u_{1})\beta_{1}+\kappa_{2}b_{4,m}=0.	
\end{equation}

Inserting all these relations into the second equation, one gets
\begin{equation*}
\begin{split}
& (g_{1}(u)v+g_{2}(u))_{t}+ (g_{3}(u)v+g_{4}u))_{x} \\
\approx & (-c[g_{1}(G)H+g_{2}(G)]+[g_{3}(G)H+g_{4}(G)]
+s'(t)(g_{1}(u_{0}\beta_{0}+g_{1}(u_{1})\beta_{1}))\delta \\
& + s(t)(g_{1}(u_{0})\beta_{0}+g_{1}(u_{1})\beta_{1}+
g_{3}(u_{0})\beta_{0}+g_{3}(u_{1})\beta_{1} \\
& +\kappa_{1}(u_{1}\alpha_{0}^{m-1}\gamma_{0}+u_{0}\alpha_{1}^{m-1}\gamma_{1})
m b_{3,m-1}+\kappa_{2}b_{4,m})\delta'\approx 0.
\end{split}
\end{equation*}
The function $s$ must be a linear one, say $s'(t)=\sigma$, and 
the above functional equation gives the last two equations in ${\mathbb R}$,
\begin{equation} \label{esp2}
-c[g_{1}(G)H+g_{2}(G)]+[g_{3}(G)H+g_{4}(G)]
+\sigma (g_{1}(u_{0})\beta_{0}+g_{1}(u_{1})\beta_{1})=0
\end{equation}
and 
\begin{equation}\label{esp3}
\begin{split}
&-c((g_{1}(u_{0})+g_{3}(u_{0})\beta_{0}+(g_{1}(u_{1})+g_{3}(u_{1})\beta_{1})\\
&+\kappa_{1}(u_{1}\alpha_{0}^{m-1}\gamma_{0}+u_{0}\alpha_{1}^{m-1}\gamma_{1})
m b_{3,m-1}+\kappa_{2}b_{4,m}=0. 
\end{split}
\end{equation}
In the above equations, only important fact about $s$ is its derivative.
Thus one can safely put $s(t)=\sigma t +\zeta$ and if the above
system (\ref{esp1}-\ref{esp3}) has a solution, then
$(u_{1},v_{1})\in \mathop{\rm SDSL}_{\zeta}(u_{0},v_{0})$ provided
that $\tilde{s}$ and $\tilde{\tilde{s}}$ can be recovered. This is certainly 
the case when $\zeta>0$. If $m$ is an odd number, then $\tilde{s}
=s(t)^{1/m}$ and $\tilde{\tilde{s}}=\tilde{s}$ are always determined.

The second part of the assertion, that $\beta_{i}$, $i=1,2$ are
independent of $\zeta$ is obvious from the above.   
\end{proof}

\begin{remark}
From the proof of the lemma one can see that it is actually possible
for $\zeta$ to take negative values, i.e.\
it is enough that $\zeta\geq -s(T)$, where $T$ is a time of interaction
when new initial data are given. 
\end{remark} 

The following assertion is crucial for the construction of weak solution 
(a solution in an associated sense)
to (\ref{gdss1}) after an interaction: At an interaction point of a DSSW
and some other wave one can consider  
the new initial value problem which contains delta function.
\medskip 

Suppose that the initial data are given by 
\begin{equation} \label{3id}
u(x,0)=\begin{cases} u_{0},& \; x<a \\ 
u_{1},& \; a<x<b \\ u_{2},& \; x>b \end{cases} \text{ and }
v(x,0)=\begin{cases} v_{0},& \; x<a \\ 
v_{1},& \; a<x<b \\ v_{2},& \; x>b \end{cases}
\end{equation}
such that there exist a singular shock wave starting from the point $x=a$
and a shock wave (or another singular shock wave) starting from the point 
$x=b$, $a<b$. They can interact if $c_{1}>c_{2}$, where $c_{i}$ is
the speed of the $i$-th wave, $i=1,2$. For the simplicity we shall assume
that $b=0$.

Let $(X,T)$ be the interaction point of the overcompressive 
singular shock wave starting at the point $x=a$
\begin{equation} \label{prvidssw}
\begin{split}
u^{1}(x,t) = & G^{1}(x-c_{1}t-a)
+\tilde{s}^{1}(t)\Big(\alpha^{1}_{0}d^{-}(x-c_{1}t-a)
+\alpha^{1}_{1}d^{+}(x-c_{1}t-a)\Big) \\
v^{1}(x,t) = & H^{1}(x-c_{1}t-a)
+ s^{1}(t)\Big(\beta^{1}_{0}D^{-}(x-c_{1}t-a)
+\beta^{1}_{1}D^{+}(x-c_{1}t-a)\Big) \\
& +\tilde{\tilde{s}}^{1}(t)\Big(\gamma^{1}_{0}d^{-}(x-c_{1}t-a)
+\gamma^{1}_{1}d^{+}(x-c_{1}t-a)\Big)
\end{split}
\end{equation}
and the admissible (singular) shock wave 
\begin{equation} \label{drugidssw}
\begin{split}
u^{2}(x,t) = & G^{2}(x-c_{2}t)
+\tilde{s}^{2}(t)\Big(\alpha^{2}_{0}d^{-}(x-c_{2}t)
+\alpha^{2}_{1}d^{+}(x-c_{2}t)\Big) \\
v^{2}(x,t) = & H^{2}(x-c_{2}t)
+ s^{2}(t)\Big(\beta^{2}_{0}D^{-}(x-c_{2}t)
+\beta^{2}_{1}D^{+}(x-c_{2}t)\Big)\\
& +\tilde{\tilde{s}}^{2}(t)\Big(\gamma^{2}_{0}d^{-}(x-c_{2}t-a)
+\gamma^{2}_{1}d^{+}(x-c_{2}t-a)\Big)
\end{split}
\end{equation}
where $G^{1}$, $G^{2}$, $H^{1}$ and $H^{2}$ are the generalized step functions 
with values $(u_{0},u_{1})$, $(u_{1},u_{2})$, $(v_{0},v_{1})$
and $(v_{1},v_{2})$, respectively. Also, $(\alpha_{0}^{i})^{m_{1}}+
(\alpha_{1}^{i})^{m_{1}}=(\gamma_{0}^{i})^{m_{1}}+
(\gamma_{1}^{i})^{m_{1}}=\beta_{0}^{i}+\beta_{1}^{i}=1$, $i=1,2$. 
Here, $m_{1}=m$ if singular part of singular shock wave is $m$SD-function
and $m_{1}=m-1$ in the case of $m'$SD-function.
If the second wave is a shock one,
then one can put $s^{2}\equiv \tilde{s}^{2}\equiv\tilde{\tilde{s}}^{2}\equiv 0$.

The speed of a singular shock wave (as well as for a shock wave) can be found
using the first equation in (\ref{gdss1}) because of assumptions (\ref{deg})
or (\ref{deg'}). For the first singular
shock wave (\ref{prvidssw}) we have 
\begin{equation*}
\begin{split}
& (f_{2}(u))_{t}+(f_{3}(u)v+f_{4}(u))_{x} \approx
(f_{2}(G))_{t}+(f_{3}(G)H+f_{4}(G))_{x} 
+ (\mathop{\rm const} s^{1}(t)\delta)_{x} \\
 \approx & (-c_{1}[f_{2}(G)]+[f_{3}(G)H+f_{4}(G)])\delta+
\mathop{\rm const} s^{1}(t)\delta'\approx 0,
\end{split}
\end{equation*}
where the term $\mathop{\rm const}s^{1}(t)$ is determined, but we shall
not write the exact value since it is not needed for the assertion.
Missing argument in the above expression is $x-c_{1}t-a$.

Let $\Gamma_{1}=\{x=c_{1}t+a\}$ and $\Gamma_{2}=\{x=c_{2}t\}$. 
Then $[\cdot]_{\Gamma_{i}}$ denotes the jump at the 
curve $\Gamma_{i}$, $i=1,2$.
Thus, one can see that the speed of that singular shock wave 
has the same value as in the case of shock wave,
$$c_{1}={ [f_{3}(G)H+f_{4}(G)]_{\Gamma_{1}}\over [f_{2}(G)]_{\Gamma_{1}}}.$$
Also, 
$$c_{2}={ [f_{3}(G)H+f_{4}(G)]_{\Gamma_{2}}\over [f_{2}(G)]_{\Gamma_{2}}}.$$

Finally, one can see that the waves given by (\ref{prvidssw}) and 
(\ref{drugidssw}) will interact at the point $(X,T)$ if $a<0$ and 
$c_{1}>c_{2}$, where 
\begin{equation*}
\begin{split}
T=& {a[f_{2}(G)]_{\Gamma_{1}}[f_{2}(G)]_{\Gamma_{2}}
\over [f_{3}(G)H+f_{4}(G)]_{\Gamma_{2}}[f_{2}(G)]_{\Gamma_{1}}
-[f_{3}(G)H+f_{4}(G)]_{\Gamma_{1}}[f_{2}(G)]_{\Gamma_{2}}}\\
X=& {a [f_{3}(G)H+f_{4}(G)]_{\Gamma_{2}}[f_{2}(G)]_{\Gamma_{1}}
\over  [f_{3}(G)H+f_{4}(G)]_{\Gamma_{2}}[f_{2}(G)]_{\Gamma_{1}}
-[f_{3}(G)H+f_{4}(G)]_{\Gamma_{1}}[f_{2}(G)]_{\Gamma_{2}}}.
\end{split}
\end{equation*}
\medskip
\medskip

Denote by $(\tilde{u}(x,t),\tilde{v}(x,t))$ a solution before interaction
time $t=T$ consisting of waves (\ref{prvidssw},\ref{drugidssw}). 

\begin{remark}
In the case of system (\ref{dss1})
one can easily calculate speeds of the above shocks and
coordinates of the interaction point.
The speeds of singular shock and entropy shock wave are
\begin{equation*}
c_{1}={u_{1}^{2}-v_{1}-u_{0}^{2}+v_{0} \over 
u_{1}-u_{0}} \text{ and } c_{2}={u_{2}^{2}-v_{2}-u_{1}^{2}+v_{1} \over 
u_{2}-u_{1}}. 
\end{equation*}
If $c_{1}>c_{2}$, then one gets 
\begin{equation*}
X={-ac_{2} \over c_{2}-c_{1}} \text{ and }
T={a \over c_{2}-c_{1}}.
\end{equation*}
for the interaction point $(X,T)$.
\end{remark}

\begin{theorem} \label{glavna}
Let system (\ref{gdss1}) be given.
Suppose that $(u_{2},v_{2})\in \mathop{\rm SDSL}_{\zeta}(u_{0},v_{0})$,
$\zeta=(\zeta_{1}+\zeta_{2})/(g_{1}(u_{0})\beta_{0}+g_{1}(u_{1})\beta_{1})$,
where the constants $\zeta_{i}$, $i=1,2$,
are defined by
\begin{equation*}
\begin{split}
& g_{1}(u^{1})v^{1}+g_{2}(u^{1})|_{(t=T)}
\approx \zeta_{1}\delta_{(X,T)} \\
& g_{1}(u^{2})v^{2}+g_{2}(u^{2})|_{(t=T)}
\approx \zeta_{2}\delta_{(X,T)}. 
\end{split}
\end{equation*}
The corresponding DSSW, $(\hat{u},\hat{v})(x,t)$ is given by 
\begin{equation} \label{izlaznidssw}
\begin{split}
\hat{u}(x,t) =  & G(x-X-c(t-T)) \\
& +\tilde{s}(t)\Big(\alpha_{0}d^{-}(x-X-c(t-T))
+\alpha_{1}d^{+}(x-X-c(t-T))\Big) \\
\hat{v}(x,t) = & H(x-X-c(t-T)) \\
& + s(t)\Big(\beta_{0}D^{-}(x-X-c(t-T))
+\beta_{1}D^{+}(x-X-c(t-T))\Big) \\
& +\tilde{\tilde{s}}(t)\Big(\gamma_{0}d^{-}(x-X-c(t-T))
+\gamma_{1}d^{+}(x-X-c(t-T))\Big) 
\end{split}
\end{equation}
for $t>T$. By Lemma \ref{podskup}, $\beta_{0}$ and $\beta_{1}$ are 
determined independently on $\zeta$, so the definition of $DSSW$ makes 
sense. 

Then there exist a solution to (\ref{gdss1},\ref{3id}) in the association 
sense such that it equals $(\tilde{u},\tilde{v})(x,t)$ 
for $t<T-\varepsilon$, and it equals $(\hat{u},\hat{v})(x,t)$
for $t>T+\varepsilon$.
\end{theorem}
\begin{proof}
Take a constant $t_{0}$ such that singular parts of the waves
$(u_{\varepsilon}^{1}(x,t),v_{\varepsilon}^{1}(x,t))$ and 
$(u_{\varepsilon}^{2}(x,t),v_{\varepsilon}^{2}(x,t))$ has disjoint
supports (i.e.\
$c_{1}t-a-c_{2}t>4\varepsilon$, for $t<T-t_{0}\varepsilon$, if
one uses the construction of the S$\delta$, $m$SD and $m'$SD-functions
defined above).

Let us denote 
\begin{equation*}
\begin{split}
& \Delta_{\varepsilon}=\{(x,t):\; |x-X|\leq t_{0}\varepsilon+\varepsilon,\;
|t-T|\leq t_{0}\varepsilon+\varepsilon\}, \\
& \tilde{\Delta}_{\varepsilon}=\{(x,t):\; |x-X|\leq t_{0}\varepsilon,\;
|t-T|\leq t_{0}\varepsilon\}, \\
& A_{\varepsilon}=\{(x,t):\; |x-X|\leq t_{0}\varepsilon+\varepsilon,\;
t=T-t_{0}\varepsilon-\varepsilon\}, \\
& B_{\varepsilon}=\{(x,t):\; x=X+t_{0}\varepsilon+\varepsilon,\;
|t-T|\leq t_{0}\varepsilon+\varepsilon\}, \\
& C_{\varepsilon}=\{(x,t):\; |x-X|\leq t_{0}\varepsilon+\varepsilon,\;
t=T+t_{0}\varepsilon+\varepsilon\}, \\
& D_{\varepsilon}=\{(x,t):\; x=X-t_{0}\varepsilon-\varepsilon,\;
|t-T|\leq t_{0}\varepsilon+\varepsilon\}.
\end{split}
\end{equation*}

Define a cut-off function $\xi_{\varepsilon}(x,t)$ which equals zero for
$(x,t)\in \Delta_{\varepsilon}$ and 1 for 
$(x,t)\in \tilde{\Delta}_{\varepsilon}$.
Let 
\begin{equation*}
(u_{temp},v_{temp})(x,t)=
\begin{cases}
(\tilde{u}(x,t),\tilde{v}(x,t)), & t<T \\
(\hat{u}(x,t),\hat{v}(x,t)), & t>T.
\end{cases}
\end{equation*}

We shall prove that the generalized functions $u$ and $v$ represented by 
\begin{equation} \label{resenje}
u_{\varepsilon}(x,t)=u_{temp}(x,t)\xi_{\varepsilon}(x,t), \text{ and }
v_{\varepsilon}(x,t)=v_{temp}(x,t)\xi_{\varepsilon}(x,t),
\; x\in {\mathbb R},\; t\geq 0
\end{equation}
solve (\ref{gdss1}) in the association sense.

Denote 
\begin{equation*}
{\mathbf F}(u,v)=\left[ \begin{matrix}f_{2}(u) \\
g_{1}(u)v+g_{2}(u) \end{matrix} \right] \text{ and }
{\mathbf G}(u,v)=\left[ \begin{matrix}f_{3}(u)v+f_{4}(u) \\
g_{3}(u)v+g_{4}(u) \end{matrix}\right].
\end{equation*} 

We have
\begin{equation*}
\begin{split}
& \iint_{{\mathbb R}_{+}^{2}}{\mathbf F}(u,v)\Psi_{t} + 
{\mathbf G}(u,v)\Psi_{x} dx dt \\
= & \iint_{\tilde{\Delta}_{\varepsilon}}{\mathbf F}(u,v)\Psi_{t} + 
{\mathbf G}(u,v)\Psi_{x} dx dt \\
= & \iint_{{\mathbb R}_{+}^{2}\setminus \tilde{\Delta}_{\varepsilon}}
{\mathbf F}(u,v)\Psi_{t} + {\mathbf G}(u,v)\Psi_{x} dx dt, 
\end{split}
\end{equation*}
for every test function $\Psi=\left[ \begin{matrix}\psi_{1}\\
\psi_{2}\end{matrix}\right]\in {\mathcal C}_{0}^\infty({\mathbb R}_{+}^{2})$.

The measure of the set $\tilde{\Delta}_{\varepsilon}$ is 
${\mathcal O}(\varepsilon^{2})$, as $\varepsilon \rightarrow 0$, 
while 
\begin{equation*}
\| {\mathbf F}(u,v)\Psi_{t}
+{\mathbf G}(u,v)\Psi_{x}\|_{L^{\infty}({\mathbb R}_{+}^{2})}
\leq \mathop{\rm const} \varepsilon^{-1+1/m}
\end{equation*}
due to the assumptions in Definition \ref{s-d}.
Thus,
\begin{equation*}
\iint_{\tilde{\Delta}_{\varepsilon}}{\mathbf F}(u,v)\Psi_{t} + 
{\mathbf G}(u,v)\Psi_{x} dx dt \sim \varepsilon^{1/m} \rightarrow 0, \text{ as }
\varepsilon \rightarrow 0. 
\end{equation*}

Using the divergence theorem for the second integral one gets
\begin{equation*}
\begin{split}
& \iint_{{\mathbb R}_{+}^{2}\setminus \tilde{\Delta}_{\varepsilon}}
{\mathbf F}(u,v)\Psi_{t} + {\mathbf G}(u,v)\Psi_{x} dx dt \\
= & \int_{\partial\tilde{\Delta}_{\varepsilon}}{\mathbf F}(u,v)\Psi\nu_{t} + 
\int {\mathbf G}(u,v)\Psi\nu_{x} ds \\
& - \iint_{{\mathbb R}_{+}^{2}\setminus \tilde{\Delta}_{\varepsilon}}
{\mathbf F}(u,v)_{t}\Psi + {\mathbf G}(u,v)_{x}\Psi dx dt. 
\end{split}
\end{equation*}
The last integral in the above expression tend to zero as 
$\varepsilon \rightarrow 0$ since $(u,v)$ solves (\ref{gdss1})
in ${\mathbb R}_{+}^{2}\setminus\tilde{\Delta}_{\varepsilon}$
due to the construction. For the other integral one gets
\begin{equation*}
\begin{split}
& \int_{\partial\tilde{\Delta}_{\varepsilon}}{\mathbf F}(u,v)\Psi \nu_{t} + 
\int {\mathbf G}(u,v)\Psi \nu_{x}ds \\
= & \int_{A_{\varepsilon}}{\mathbf F}(u,v)\Psi dx 
- \int_{C_{\varepsilon}}{\mathbf F}(u,v)\Psi dx
+ \int_{D_{\varepsilon}}{\mathbf G}(u,v)\Psi dt 
- \int_{B_{\varepsilon}}{\mathbf G}(u,v)\Psi dt.
\end{split}
\end{equation*}
Functions $u_{\varepsilon}$ and $v_{\varepsilon}$ are $L^{\infty}$-bounded
uniformly in $\varepsilon$ on the sides 
$B_{\varepsilon}$ and $D_{\varepsilon}$. Since their lengths
are ${\mathcal O}(\varepsilon)$, integrals over them tends to zero
as $\varepsilon \rightarrow 0$. 

Using the fact that $f_{2}(d_{\varepsilon})\approx 0$ one gets
\begin{equation*}
\lim_{\varepsilon \rightarrow 0}F(\tilde{u},\tilde{v})|_{t=T}
=\left[\begin{matrix} 0 \\
(\zeta_{1}+\zeta_{2})\delta_{(X,T)}\end{matrix}\right],
\end{equation*} 
as well as the construction of S$\delta$- and $m'$SD (or $m$SD)-functions,
one gets
\begin{equation*}
\lim_{\varepsilon \rightarrow 0}
\int_{A_{\varepsilon}} {\mathbf F}(u_{\varepsilon},v_{\varepsilon})dx
=\left[\begin{matrix} 0 \\
\zeta_{1}+\zeta_{2}\end{matrix}\right]\cdot \Psi(X,T).
\end{equation*} 
Thus, there has to be true that
\begin{equation*}
\lim_{\varepsilon \rightarrow 0}
\int_{C_{\varepsilon}} {\mathbf F}(u_{\varepsilon},v_{\varepsilon})dx
=-\left[\begin{matrix} 0 \\
\zeta_{1}+\zeta_{2}\end{matrix}\right]\cdot \Psi(X,T).
\end{equation*} 
This implies $f_{2}(\hat{u})|_{(X,T)}\approx 0$ and  
\begin{equation} \label{dusl}
g_{1}(\hat{u})\hat{v}+g_{2}(\hat{u})|_{(X,T)}
\approx (\zeta_{1}+\zeta_{2})\delta_{(X,T)}.
\end{equation}
Due to conditions (\ref{deg}) or (\ref{deg'}) one immediately
gets $f_{2}(\hat{u})|_{(X,T)}\approx 0$. Put 
$\zeta=(\zeta_{1}+\zeta_{2})/(g_{1}(u_{0})\beta_{0}+g_{1}(u_{1})\beta_{1})$.
Then
$$
g_{1}(\hat{u})\hat{v}+g_{2}(\hat{u})|_{t=T}
\approx \hat{G}\hat{H}+\hat{G}
+s(T)(g_{1}(u_{0})\beta_{0}+g_{1}(u_{1})\beta_{1})\delta(X)
$$
and after another restriction on $x=X$,
$$
g_{1}(\hat{u})\hat{v}+g_{2}(\hat{u})|_{(X,T)}\approx 
(\zeta_{1}+\zeta_{2})\delta_{(X,T)}.
$$
This concludes the proof.
\end{proof}

\begin{remark} \label{distriblim}
The distributional limit of the result of the interaction is given by
\begin{equation*}
\begin{split}
u(x,t) 
&= \left\{ \begin{aligned} 
u_{0},&\; x<c_{1}t-a,\;t<t \\
u_{1},&\; c_{1}t-a<x<c_{2}t, \; t<T \\
u_{2},&\; x>c_{2}t, \; t<T \\
u_{0},&\; x<ct+X, \; t>T \\
u_{2},&\; x>ct+X, \; t>T 
\end{aligned} \right.  \\ 
v(x,t) 
&= \left\{ \begin{aligned} 
v_{0},&\; x<c_{1}t-a,\;t<t \\
v_{1},&\; c_{1}t-a<x<c_{2}t, \; t<T \\
v_{2},&\; x>c_{2}t, \; t<T \\
v_{0},&\; x<ct+X, \; t>T \\
v_{2},&\; x>ct+X, \; t>T 
\end{aligned} \right\} + s_{1}(t)\delta_{S_{1}}+
s_{2}(t)\delta_{S_{2}}+s(t)\delta_{S},
\end{split}
\end{equation*}
where $S_{1}=\{(x,t):\; x=c_{1}t+a,\; t\in [0,T]$, 
$S_{2}=\{(x,t):\; x=c_{2}t,\; t\in [0,T]$ and 
$S=\{(x,t):\; x-X=c(t-T),\; t\in [T,\infty)$. 
If the second wave (\ref{drugidssw}) is a shock one, then $s_{2}\equiv 0$. 

The above solution is continuous in
$t$ with values in ${\mathcal D}'({\mathbb R})$. This fact 
can be used in the approach similar to \cite{ShDan1}, where 
the variable $t$ is treated separately, i.e.\
when system (\ref{gdss1}) is considered to be in evolution form.
\end{remark}

The theorem shows that after an interaction of a singular shock with 
some shock or another singular shock the problem 
reduces to solving system (\ref{gdss1}) with the new initial data
(\ref{eq6}). 
\begin{remark}

\noindent 
(i) The solution to the interaction problem from Theorem \ref{glavna} 
is always associated with a lower association rate 
than the solution of the original Riemann problem.
For specific system it seems possible to make more sophisticated
construction in order to improve the rate.

\noindent
(ii) It appears that $d_{\varepsilon}^{\pm}$ 
are unavoidable  correction factors even their distributional limit 
equals zero.

The conditions (\ref{deg}) and (\ref{deg'}) ensures that the new initial data
at intersection point do not depend on $m$SD- or $m'$SD-functions in 
the solution. We have used them because the real nature of 
$m$SD- and $m'$SD-functions is not so clear yet.
\end{remark}

The above theorem will be used in 
the rest of the paper for investigation of interactions between singular
shock waves and other types of waves in the special case 
of system (\ref{dss1}).

\section{Applications}
Consider now system (\ref{dss1}) which a special case to (\ref{gdss1}).
The authors of \cite{KeyKr} defined and proved existence 
of singular shock wave solutions for some Riemann problems 
of this system. 

In the present paper, we will investigate interactions of 
such solutions with the other solutions to Riemann problem for (\ref{dss1}).
In order to familiarize a reader with the presented results, let us 
give some basic remarks about such solutions.

For a given Riemann data $(u_{0},v_{0})$, $(v_{0},v_{1})$, there are three
basic solution types:
\begin{enumerate}
\item[(a)] {\it Shock waves}
\begin{equation}\label{sw12}
u(x,y)=\left\{ \begin{aligned} u_{0},& \; x<ct \\ u_{1},& \;
x>ct \end{aligned} \right. \phantom{second}
v(x,y)=\left\{ \begin{aligned} v_{0},& \; x<ct \\ v_{1},& \;
x>ct \end{aligned} \right. 
\end{equation}
where $c=[u^{2}-v]/[u]$ and $(u_{1},v_{1})$ lies in an admissible
part of Hugoniot locus of the point $(u_{0},v_{0})$.
\item[(b)] {\it Centered rarefaction waves}
\begin{equation}\label{rw1}
\begin{split}
& u(x,t)=\left\{ \begin{aligned} u_{0},& \; x<(u_{0}-1)t \\
x/t+1,& \; (u_{0}-1)t\leq x \leq (u_{1}-1)t \\
u_{1},& \; x>(u_{1}-1)t
\end{aligned} \right. \\
& v(x,t)=\left\{ \begin{aligned} v_{0},& \; x<(u_{0}-1)t \\
(x/t)^{2}/2+2x/t+C_{1},& \; (u_{0}-1)t\leq x \leq (u_{1}-1)t \\
v_{1},& \; x>(u_{1}-1)t
\end{aligned}\right.
\end{split}
\end{equation}
(1-rarefaction wave), where $C_{1}=v_{0}-u_{0}^{2}/2-u_{0}-1/2$, 
when $(u_{1},v_{1})$ lies in an 1-rarefaction curve
starting at the point $(u_{0},v_{0})$. Or 
\begin{equation}\label{rw2}
\begin{split}
& u(x,t)=\left\{ \begin{aligned} u_{0},& \; x<(u_{0}+1)t \\
x/t-1,& \; (u_{0}+1)t\leq x \leq (u_{1}+1)t \\
u_{1},& \; x>(u_{1}+1)t
\end{aligned}\right. \\
& v(x,t)=\left\{ \begin{aligned} v_{0},& \; x<(u_{0}+1)t \\
(x/t)^{2}/2-2x/t+C_{2},& \; (u_{0}+1)t\leq x \leq (u_{1}+1)t \\
v_{1},& \; x>(u_{1}+1)t
\end{aligned}\right.
\end{split}
\end{equation}
(2-rarefaction wave), where $C_{2}=v_{0}-u_{0}^{2}/2+u_{0}-1/2$, 
when $(u_{1},v_{1})$ lies in an 2-rarefaction curve
starting at the point $(u_{0},v_{0})$.
\item[(c)] {\it Singular shock waves} (see Definition 
\ref{prom}) of $3'$SD-type,
\begin{equation}\label{ssw}
\begin{split}
& u(x,y)=\left\{ \begin{aligned} u_{0},& \; x<ct \\ u_{1},& \;
x>ct \end{aligned} \right\}+\tilde{s}(t)(\alpha_{0}d_{\varepsilon}^{-}(x-ct)
+ \alpha_{1}d_{\varepsilon}^{+}(x-ct))\\
& v(x,y)=\left\{ \begin{aligned} v_{0},& \; x<ct \\ v_{1},& \;
x>ct \end{aligned} \right\}+s(t)(\beta_{0}D_{\varepsilon}^{-}(x-ct)+
\beta_{1}D_{\varepsilon}^{+}(x-ct)),
\end{split} 
\end{equation}
where $c=[u^{2}-v]/[u]$, and all other terms 
are as in that definition. That means
\begin{equation}\label{usl}
D_{\varepsilon}\approx \delta,
\;(d_{\varepsilon}^{\pm})^{i}\approx 0, \; i=1,3, \; (d_{\varepsilon}^{\pm})^{2}
\approx \delta,
\end{equation} 
while $(u_{1},v_{1})$ lies in a region denoted by $Q_{7}$ in \cite{KeyKr}
of the point $(u_{0},v_{0})$ (see Figure 1). 
\end{enumerate}
For an arbitrary Riemann problem to (\ref{dss1})
one can construct a solutions by 
the means of these waves or their combinations (\cite{KeyKr}).
\medskip 

While interactions of the first two types can be handled in a usual
way, interactions involving singular shock waves are quite different 
and far more interesting, so they become a topic of this paper.

The procedure for the singular shock wave interactions
can be also used for systems (\ref{gdss1}). But
a complete after-interaction solution highly depends on a particular 
system. That is the reason why we treat system (\ref{dss1}) only.

In order to simplify notation, we shall substitute the point  
$(X,T)$ in (\ref{eq6}) by $(0,0)$ and then solve the
Cauchy problem (\ref{dss1},\ref{eq6}).

There are no multiplication of $v$ with $u$ 
in system (\ref{dss1}), so in the sequel it will be enough 
to take $D^{-}=D^{+}$, 
$\alpha_{0}(t):=\alpha_{0}\tilde{s}(t)$, 
$\alpha_{1}(t):=\alpha_{1}\tilde{s}(t)$ and
$\beta(t):=s(t)$, i.e.\ 
to look for a solution of the form
\begin{equation} \label{eq7}
\begin{split}
& u = G(x-ct)+(\alpha_{0}(t)d^{-}(x-ct)
+\alpha_{1}(t)d^{+}(x-ct)) \\
& v = H(x-ct)+\beta(t)D(x-ct),
\end{split}
\end{equation}
where $G$ and $H$ are generalized step functions, 
while $d$ is $3'$SD- and $D$ is S$\delta$-function and $c\in {\mathbb R}$.

Let us determine SDSL of (\ref{dss1}) for some $(u_{0},v_{0})\in 
{\mathbb R}^{2}$.  

Substitution of (\ref{eq7}) into the first equation of the system gives
\begin{equation}\label{eq8} 
\begin{split}
& c={u_{1}^{2}-v_{1}-u_{0}^{2}+v_{0} \over u_{1} - u_{0}} \\
& \alpha_{0}^{2}(t)+\alpha_{1}^{2}(t)=\beta(t),
\end{split}
\end{equation}
where $c$ is the speed of the wave. After neglecting all terms 
converging to zero as $\varepsilon \to 0$, 
the second equation becomes
\begin{equation*}
\begin{split}
& \partial_{t}H_{\varepsilon}(x-ct)+\beta'(t)\delta(x-ct)
-c\beta(t)\delta'(x-ct)
+\partial_{x}({1 \over 3}G_{\varepsilon}^{3}-G) \\ 
&+(u_{1}\alpha_{0}^{2}(t)+u_{0}\alpha_{1}^{2}(t))\delta'(x-ct)= 0. 
\end{split}
\end{equation*}
Thus, the following relations has to hold.
\begin{equation} \label{eq8dva}
\beta'(t)=c(v_{1}-v_{0})-\big( {1 \over 3} u_{0}^{3} - u_{0}
-{1 \over 3} u_{1}^{3} + u_{1} \big) =: k,
\end{equation}
i.e.
\begin{equation*}
\beta(t)=kt+\zeta, \mbox{ since } \beta(0)=\zeta
\end{equation*}
and 
\begin{equation} \label{eq9}
u_{1}\alpha_{0}^{2}(t)+u_{0}\alpha_{1}^{2}(t)=c\beta(t).
\end{equation}

Like in \cite{KeyKr} one can see that the overcompressibility means
\begin{equation*}
u_{0}-1 \geq c \geq u_{1}+1,
\end{equation*}
i.e., $v_{1}$ lies between the curves
\begin{equation*}
\begin{split}
& D=\{ (u,v):\; v=v_{0}+u^{2}+u-u_{0}u-u_{0} \} \\
& E=\{ (u,v):\; v=v_{0}-u+u_{0}u-u_{0}^{2}+u_{0} \}, 
\end{split}
\end{equation*}
and $u_{0}-u_{1}\geq 2$.

\begin{figure}[Ht]
\includegraphics*[scale=0.75]{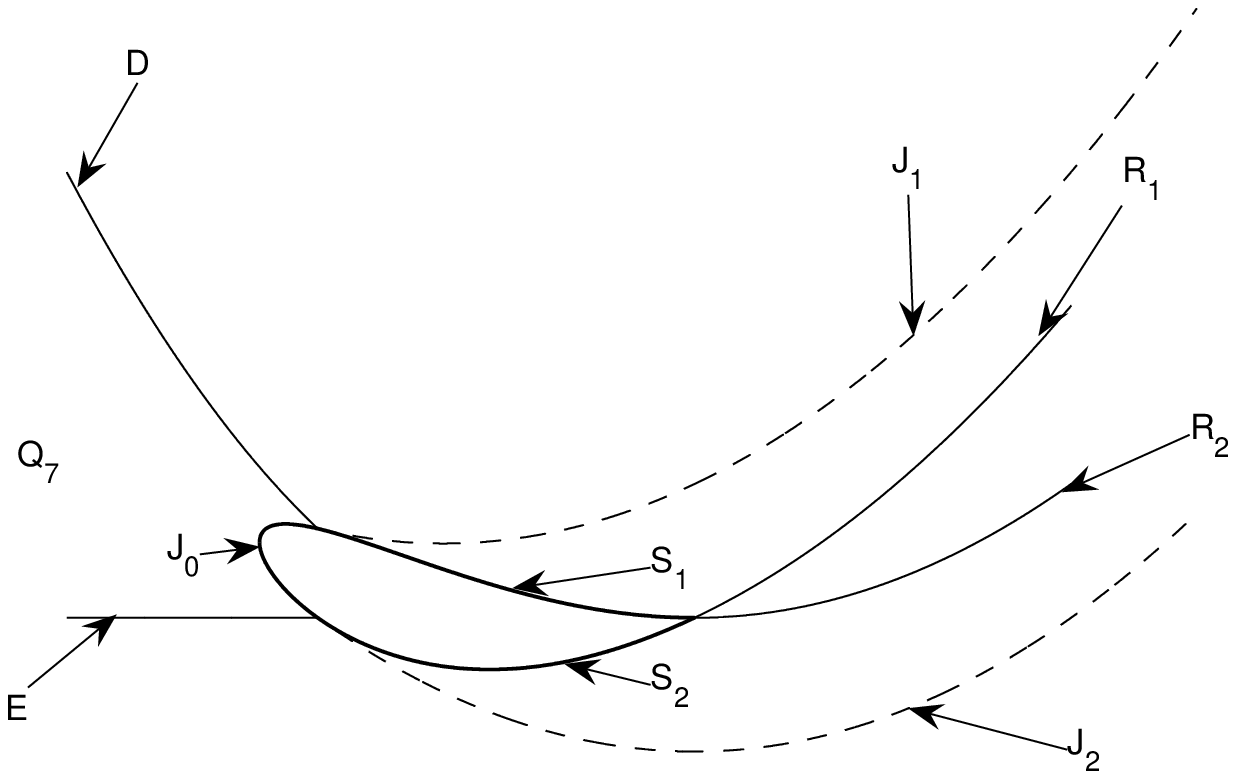} 
%
%
\caption{}
\end{figure} 

Denote by $J_{1}$ the union of the parts of admissible 
Hugoniot locus 
\begin{equation*}
S_{1}=\Big \{ (u,v): \;
v-v_{0}=(u-u_{0})\Big( {u_{0}+u \over 2}
+\sqrt{1-{(u_{0}-u)^{2}\over 12}}\Big)\Big\}, 
\end{equation*}
and
\begin{equation*}
S_{2}=\{ (u_{1},v_{1}): \;
v-v_{0}=(u-u_{0})\Big( {u_{0}+u \over 2}
-\sqrt{1-{(u_{0}-u)^{2}\over 12}}\Big)\Big\},
\end{equation*}
for $u\in [u_{0}-\sqrt{12},u_{0}-3]$. Note that $S_{i}$ is not an
$i$th shock curve but only a label. 

The points between the curves 
$D$ and $E$, and on the left-hand side
of $J_{1}$ defines the area denoted by $Q_{7}$ in \cite{KeyKr}. 
Here, this area is called delta singular locus.  

One can easily check that system
(\ref{eq8},\ref{eq9}) has a solution if and only if $\beta(t)>0$.
\medskip 
Depending on $k$, defined in (\ref{eq8dva}), 
there are three possibilities for a resulting wave:

\noindent 
(i) If $k>0$, then $\tilde{\beta}'(t)>0$ and $(u_{1},v_{1})\in Q_{7}$. 
The resulting singular shock has the same properties as
before, i.e.\ 
its strength increases with the time.
\medskip

\noindent
(ii) If $k=0$, then $\tilde{\beta} \equiv \text{const} = \zeta >0$ 
and the corresponding part of a singular overcompressive locus is $J_{1}$.
The result of the interaction is a new kind of singular shock wave,  
its strength is a constant with respect to the time.  
\medskip 

\noindent
(iii) If $k<0$ (this means that the point $(u_{1},v_{1})$ 
is on the left-hand side of $J_{1}$), then 
the resulting singular shock wave has much more 
differences from the usual one (with an increasing strength). 
Its initial strength equals $\zeta$,
$\beta(0)=\zeta>0$, but linearly decreases in time. At some point
$T_{0}$ the strength of
the singular shock equals zero and the singular 
shock wave does not exist after that.
In the rest of the paper we shall see some cases when this happens.
The new initial data for time $t=T_{0}$ are the Riemann ones, and the 
solution after that time can be find in the usual way,
by using the results in \cite{KeyKr}.
\medskip

All the above facts are collected in the following theorem.

\begin{theorem} \label{2sl}
The SDSL$_{\zeta}$, $\zeta>0$, for (\ref{dss1},\ref{eq6}) is the 
area bounded by the curves $D$, $E$, $S_{2}\setminus J_{1}$ and   
$S_{1}\setminus J_{1}$. (The area $Q_{7}$ is a subset of this one,
as known from Lemma \ref{podskup}.)
The overcompressive SDSL$_{\zeta}$, $\zeta >0$, 
is a part of the SDSL bounded by the curves $D$ and $E$ such that
$u_{1}\leq u_{0}-2$.
\end{theorem} 

\subsection{Interaction of a singular shock and an admissible shock wave}

Suppose that a singular shock wave with a speed $c_{1}$ and 
a left- and right-hand values 
$U_{0}=(u_{0},v_{0})$ and $U_{1}=(u_{1},v_{1})$, respectively,
interact with an admissible shock wave with a speed $c_{2}<c_{1}$ 
having left-hand and right-hand values 
$U_{1}=(u_{1},v_{1})$ and $U_{2}=(u_{2},v_{2})$, respectively, 
at a point 
$(X,T)$. 
\begin{lemma} \label{dss_sw}
If the above singular shock and shock wave are admissible,
$(u_{2},v_{2})$ lies between the lines $D$ and $E$.
Thus, the solution after the interaction is a single 
overcompressive singular shock wave.
\end{lemma}
\begin{proof}
Since $u_{0} \geq u_{1}+3$ and $u_{1}>u_{2}$ 
(because of the admissibility conditions for singular and shock wave), we have
$u_{0} > u_{2}+3$. The point $(u_{2},v_{2})$ lies on the curve 
$S_{1}$ or $S_{2}$ with the origin at the point $(u_{1},v_{1})$. 
Thus
\begin{equation*}
v_{2}=v_{1}+(u_{2}-u_{1})\Big({u_{1}+u_{2}\over 2} \pm
\sqrt{1-{(u_{1}-u_{2})^{2} \over 12}}\Big).
\end{equation*}
The point $(u_{1},v_{1})$ lies in the area denoted by $Q_{7}$, 
thus bellow or at the curve $D$ with the origin at 
$(u_{0},v_{0})$. Therefore 
\begin{equation*}
v_{1}\leq v_{0}+u_{1}^{2}+u_{1}-u_{0}u_{1}-u_{0}.
\end{equation*}
Let the point $(u_{0},v_{0})$ be the origin. The point $(u_{2},v_{2})$
will be bellow the curve $D$ if 
\begin{equation*}
\begin{split}
& v_{0}+u_{1}^{2}+u_{1}-u_{0}u_{1}-u_{0}+
(u_{2}-u_{1})\Big({u_{1}+u_{2}\over 2} \pm
\sqrt{1-{(u_{1}-u_{2})^{2} \over 12}}\Big) \\
\leq & v_{0}+u_{2}^{2}+u_{2}-u_{0}u_{2}-u_{0}.
\end{split}
\end{equation*}
Non-positivity of $u_{1}-u_{2}$ gives
\begin{equation*}
\pm \sqrt{1-{(u_{1}-u_{2})^{2} \over 12}} \leq 
{1 \over 2}(u_{0}-u_{1})+{1 \over 2}(u_{0}-u_{2})-1.
\end{equation*}
The left-hand side of the above inequality is less than 2, while 
the right-hand side is greater that 2. Thus, the point $(u_{2},v_{2})$
really lies bellow the curve $D$. 

In the same way one can prove that the point $(u_{2},v_{2})$
lies above the curve $E$.
\end{proof}

\begin{remark} \label{another}
In the same manner as above, one can prove that the situation is the same
when singular shock and shock wave change sides. That is, when an admissible
singular shock wave interacts with an admissible shock wave 
from the right-hand side, then the solution is again a single admissible 
singular shock wave.
\end{remark}

\subsection{Double singular shock wave interaction}

Suppose that an admissible singular shock wave with a speed $c_{1}$ and 
left- and right-hand 
side values $U_{0}=(u_{0},v_{0})$ and $U_{1}=(u_{1},v_{1})$,
respectively, interacts 
with an another singular shock wave with a speed $c_{2}<c_{1}$ and 
left-hand (right-hand) side values
$U_{1}=(u_{1},v_{1})$ ($U_{2}=(u_{2},v_{2})$) at the point $(X,T)$.
Since the conditions for the existence of singular shock waves include
$u_{0}-u_{1}\geq 3$ and $u_{1}-u_{2}\geq 3$, then $u_{0}-u_{2}\geq 6$, i.e.\
the point $(u_{2},v_{2})$ is on the left-hand side of the line 
$u=u_{0}-\sqrt{12}$.
Concerning the position of the point $(u_{2},v_{2})$ 
in the plane of wave regions with the origin at $(u_{0},v_{0})$ there 
are three possibilities:

\noindent
(i) The point $(u_{2},v_{2})$ is between or at the curves $D$ and $E$.
The result of the interaction is a single singular shock wave (with 
increasing strength).

\noindent
(ii) The point $(u_{2},v_{2})$ is above the curve $D$.
The result of the interaction is an 1-rarefaction wave followed 
with a singular shock wave.

\noindent
(iii) The point $(u_{2},v_{2})$ is bellow the curve $E$.
The result of the interaction is a singular shock wave 
followed by a 2-rarefaction wave.

SDSL's always have increasing strength in these three cases.

\section{Intersection of a singular shock wave and a rarefaction wave}

The last possibility of singular shock wave interactions is between
a singular shock wave and a rarefaction wave.  
That possibility is omitted from a considerations of the general case
due to a richness of possible behaviors. Nevertheless, the most
of specific Riemann problems can be treated similarly as system (\ref{dss1})
was here, at least up to some point.

For a given point $(u_{0},v_{0})$, the rarefaction curves are given by 
(see \cite{KeyKr})
\begin{equation*}
\begin{split}
& R_{1}=\{(u,v):\;
v=v_{0}-{1\over 2}u_{0}^{2}+{1\over 2}u^{2}+u-u_{0}\}. \\
& R_{2}=\{(u,v):\;
v=v_{0}-{1\over 2}u_{0}^{2}+{1\over 2}u^{2}-u+u_{0}\} 
\end{split}
\end{equation*}

Suppose that a singular shock wave with left- and right-hand side values
$U_{0}=(u_{0},v_{0})$ and $U_{1}=(u_{1},v_{1})$, 
from the left-hand side interacts
with a rarefaction wave at some point $(X,T)$.
If the rarefaction wave is approximated with a number of small 
amplitude (non-admissible) shock waves like in wave fronth tracking 
algorithm (see \cite{Bre} for example),
intuition given in Theorem \ref{glavna},
such that the first task should be to look at the singular shock
wave and the interaction of singular shock and non-admissible shock wave.
It is possible to extend Theorem \ref{glavna} for such a case,
providing that a non-admissible shock wave has amplitude small enough
(of the rate $\varepsilon^{2}$, say).
Denote by $(u_{r},v_{r})$ the end-point in a rarefaction curve. Let us note
that the starting point of the curve $(u_{1},v_{1})$ is in  $Q_{7}$.

In what follows, we shall abuse the notation and denote by 
$(u_{1},v_{1})$ the left-hand side of an approximated non-admissible 
shock wave. Denote by $(u_{1},v_{1})\in Q_{7}$ the left-hand side and by 
$(u_{2},v_{2})$ the right-hand side value of a part from the rarefaction
curve. If $(u_{2},v_{2})\in Q_{7}$, then the result of the interaction is a
single singular shock wave, with the left-hand side value equals 
$(u_{0},v_{0})$. The speed depends on initial values as in (\ref{eq8}).
So, one can continue the procedure taking approximate 
points from the rarefaction 
curve as the right-hand values of the non-admissible shock wave until it
reaches the border of $Q_{7}$. 

After looking at the above discrete model we are back in a real situation.

Let us denote by $(c(t),t)$, 
$t$ belonging to  some interval, 
a path of the resulting singular shock wave trough $Q_{7}$.
It is possible to explicitly calculate
the above path. For example if a singular shock wave interacts with 
a centered 1-rarefaction waves,  
substituting 
\begin{equation*}
\begin{split}
& u(x,t)=\left\{ \begin{aligned} u_{1}, & \; x<c(t) \\
\phi_{1}(x/t), & \; x>c(t) \end{aligned} \right\} + 
\alpha_{0}(t)d_{\varepsilon}^{-}(x-c(t))
+\alpha_{1}(t)d_{\varepsilon}^{+}(x-c(t)) \\
& v(x,t)=\left\{ \begin{aligned} v_{1}, & \; x<c(t) \\
\phi_{2}(x/t), & \; x>c(t) \end{aligned} \right\} + 
\beta(t)D_{\varepsilon}(x-c(t))
\end{split}
\end{equation*}
in system (\ref{dss1}), one obtains  
\begin{equation*}
\begin{split}
& \tilde{\alpha_{0}}^{2}(t)+\tilde{\alpha}_{1}^{2}(t)=\tilde{\beta}(t)\\
& c(t)=\Big(t(1-2(u_{1}-v_{0}+v_{1}+u_{0}^{2}-u_{1}^{2}))\\
& + T(1-2(u_{0}-v_{1}-u_{0}u_{1}+u_{0}^{2}-u_{1}^{2}))\Big)/(2(u_{0}-1))\\
& \tilde{\beta}'(t)=c'(t)\Big({1\over 2}\Big({c(t) \over t}+1\Big) 
 +\Big({c(t) \over t}+1\Big)+v_{1}-{1\over 2}u_{1}^{2}-u_{1}-v_{0}\Big)\\
& -\Big({1\over 3}\Big({c(t) \over t}+1\Big)^{3}-\Big({c(t) \over t}+1\Big)
-{1\over 3}u_{0}^{3}+u_{0}\Big),
\end{split}
\end{equation*}
where the initial data for $\beta$ at the point $t=T$ is the initial strength 
of the singular shock wave $\beta(T)$.
The above calculations means that a form of the resulting singular shock
curve and its strength are uniquely determined trough the area $Q_{7}$.
If $(u_{r},v_{r})\in Q_{7}$, then the analysis is finished. Suppose that
this is not true. The main problem is to analyse situation when
rarefaction curve intersects the boundary of $Q_{7}$. Let us try to
find out what is happening by using
a discrete model.

Thus, the first real problem is to find a form of solution
when the points from the rarefaction curve satisfy:
$(u_{1},v_{1})\in Q_{7}$ and $(u_{2},v_{2})\not\in Q_{7}$.

Denote by $\tilde{D}$ and $\tilde{G}$ the intersection points of the 
curve $J_{1}$ (or the line $u=u_{0}-3$) with the curves 
$E$ and $D$, respectively (see Figure 2).

\subsection{The first critical case}

Denote by $J$ the 1-rarefaction curve starting from the point $\tilde{G}$ and
by $J_{2}$ the 2-rarefaction curve starting from the point $\tilde{D}$
The region where 
$(u_{2},v_{2})$ can lie consist of five subregions:

\noindent
(i) {\it The rarefaction curve which starts at $(u_{1},v_{1})$
intersects the curve $D$ out of point $\tilde{G}$.}
The point $(u_{2},v_{2})$ lies in the region above the curve $D$ and 
left of the line $u=u_{0}-3$. The final result of the interaction 
is a 1-rarefaction wave ($R_{1}$)
followed by a singular shock wave with increasing strength.

\noindent
(ii) {\it The rarefaction curve which starts at $(u_{1},v_{1})$
intersects the curve $E$ out of point $\tilde{D}$.}
The point $(u_{2},v_{2})$ lies in the region below the curve $E$ and 
on the left-hand side of the line $u=u_{0}-3$. The result of 
the interaction is a singular shock wave 
with increasing strength followed by a 2-rarefaction wave ($R_{2}$).

\noindent
(iii) {\it The rarefaction curve which starts at $(u_{1},v_{1})$
intersects the curve $J_{1}$ out of points $\tilde{D}$ 
and $\tilde{G}$.} Since an amplitude of a non-admissible 
shock wave can be as small as necessary, one can assume that 
the point $(u_{2},v_{2})$ lies in the 
second delta singular locus and the resulting singular shock wave
has a negative strength.
The strength-function $\tilde{\beta}(t)=\zeta+k(t-T_{0})$ of the resulting 
singular shock is decreasing, so, there could exists 
a point $T_{1}=T-\zeta/k$ such that
$\tilde{\beta}(T_{1})=\alpha_{0}=\alpha_{1}=0$. Let 
$X_{1}=cT_{1}+(X-T)$, where $c$ is the speed of the resulting 
singular shock wave (space coordinate of the point where strength
reaches zero). Therefore, in the time $t=T_{1}$, we have to solve 
new Riemann problem
\begin{equation*}
u|_{t=T_{1}}=\begin{cases} u_{0}, & x<X_{1} \\
u_{2},& x>X_{1} \end{cases}, \;
v|_{t=T_{1}}=\begin{cases} v_{0}, & x<X_{1} \\
v_{2},& x>X_{1} \end{cases}.
\end{equation*} 
This problem has a unique entropy solution consists from 
two shock waves, since the point $(u_{2},v_{2})$ is between the curves
$S_{1}$ and $S_{2}$, with respect to the origin at the 
point $(u_{0},v_{0})$. This means that the singular 
shock wave decouples into a pair of admissible shock waves.
If $u_{r}\leq u_{0}-2$, this pair of the shock waves are the final solution. 
The case when $u_{r}<u_{0}-2$ belongs to the following
subsection, i.e.\
the second critical case. 

%

\noindent
(iv) {\it The rarefaction curve $R_{j}$, $j=1 \text{ or } 2$,
which starts at $(u_{1},v_{1})$
intersects the curve $J_{1}$ in the point $\tilde{G}$.}
We can take $\tilde{G}=(u_{2},v_{2})$ for convenience. The set 
of such points  $(u_{1},v_{1})$ lies on the inverse rarefaction curve, 
which starts from the right-hand side values, i.e.\ 
\begin{equation*}
\begin{split}
&\tilde{R}_{1}=\{(u,v):\; v=v_{0}+(u_{0}^{2}-u^{2})/2+u_{0}-u\} \\
&\text{and} \\
&\tilde{R}_{2}=\{(u,v):\; v=v_{0}+(u_{0}^{2}-u^{2})/2-u_{0}+u\}
\end{split}
\end{equation*}
(the same explanation will be used in Remark \ref{r4} bellow). 
Straightforward calculation
shows that this curve lies in the region $Q_{7}$, thus this 
situation is possible, as one can see using the 
inverse rarefaction curves $\tilde{R}_{1}$ and $\tilde{R}_{2}$ 
given above. 

If $j=1$, then 
the point $(u_{2},v_{2})$ belongs
to $J$ and the solution after the interaction is an $R_{1}$-wave
followed by a singular shock wave with a constant strength.

If $j=2$, then
the point $(u_{2},v_{2})$ lies in the area bellow the curve $J$.
This can be verified by direct calculation,
taking into account that the amplitude of 
a non-admissible shock is small enough,
$u_{2}<u_{0}-2$. The solution after the interaction is an 
admissible singular shock wave with a decreasing strength.
Further explanations of a such singular shock wave is given in the
following subsection.
 
\noindent
(v) {\it The rarefaction curve $R_{j}$ which starts at $(u_{1},v_{1})$
intersects the curve $J_{1}$ in the point $\tilde{D}$.}
Again, let $\tilde{D}=(u_{2},v_{2})$.
Simple calculation, as in 
the case (iv), shows that this situation is also possible 
since the inverse rarefaction curves starting from 
$\tilde{D}$ stay in $Q_{7}$.
If $j=2$, the point $(u_{2},v_{2})$ belongs
to $J_{2}$, and then the solution after the interaction is 
singular shock wave with a constant strength followed by an
$R_{2}$-wave. 

If $j=1$, then use of the same arguments as above gives that 
the point $(u_{2},v_{2})$ lies in the area above the curve $J_{2}$
and the result of the interaction is an admissible  
singular shock wave with a decreasing strength. Again, one can see 
the following subsection for the further analysis.

\subsection{The second critical case}

Now we are dealing with the problem when the rarefaction wave
after passing through $J_{1}$, after 
passes trough the curves $D$ or $E$. 
One can see that this is the continuation of the cases (iii)-(v)
from the previous part.

\noindent
(a) Denote by $\hat{D}$ the area above the curve $D$,
bellow $S_{1}$ and on the left-hand side of the line $u=u_{0}-2$. 
Also denote by $\hat{\hat{D}}$ the area above the curve $E$,
bellow $S_{1}$ and on the right-hand side of the line $u=u_{0}-2$.
If $(u_{2},v_{2})$ lies in one of these regions, 
the solution is combination of a 
rarefaction wave $R_{1}$ and an overcompressive singular shock wave with 
a decreasing or constant strength. 

\noindent
(b) Denote by $\hat{E}$ the area bellow the curve $E$,
above $S_{2}$ and on the left-hand side of the line $u=u_{0}-2$. 
Also denote by $\hat{\hat{E}}$ the area bellow the curve $D$,
above $S_{2}$ and on the right-hand side of the line $u=u_{0}-2$.
If $(u_{2},v_{2})$ lies in one of these regions, 
the solution is then a combination of
an overcompressive singular shock wave with 
a decreasing or constant strength and a rarefaction wave $R_{2}$.
Denote by $D_{0}$ the area bounded by the curves $D$, $E$, $S_{1}$
and $S_{2}$ such that $u<u_{0}-2$ in $D_{0}$.
\begin{figure}[Ht]
%
\includegraphics*[scale=0.75]{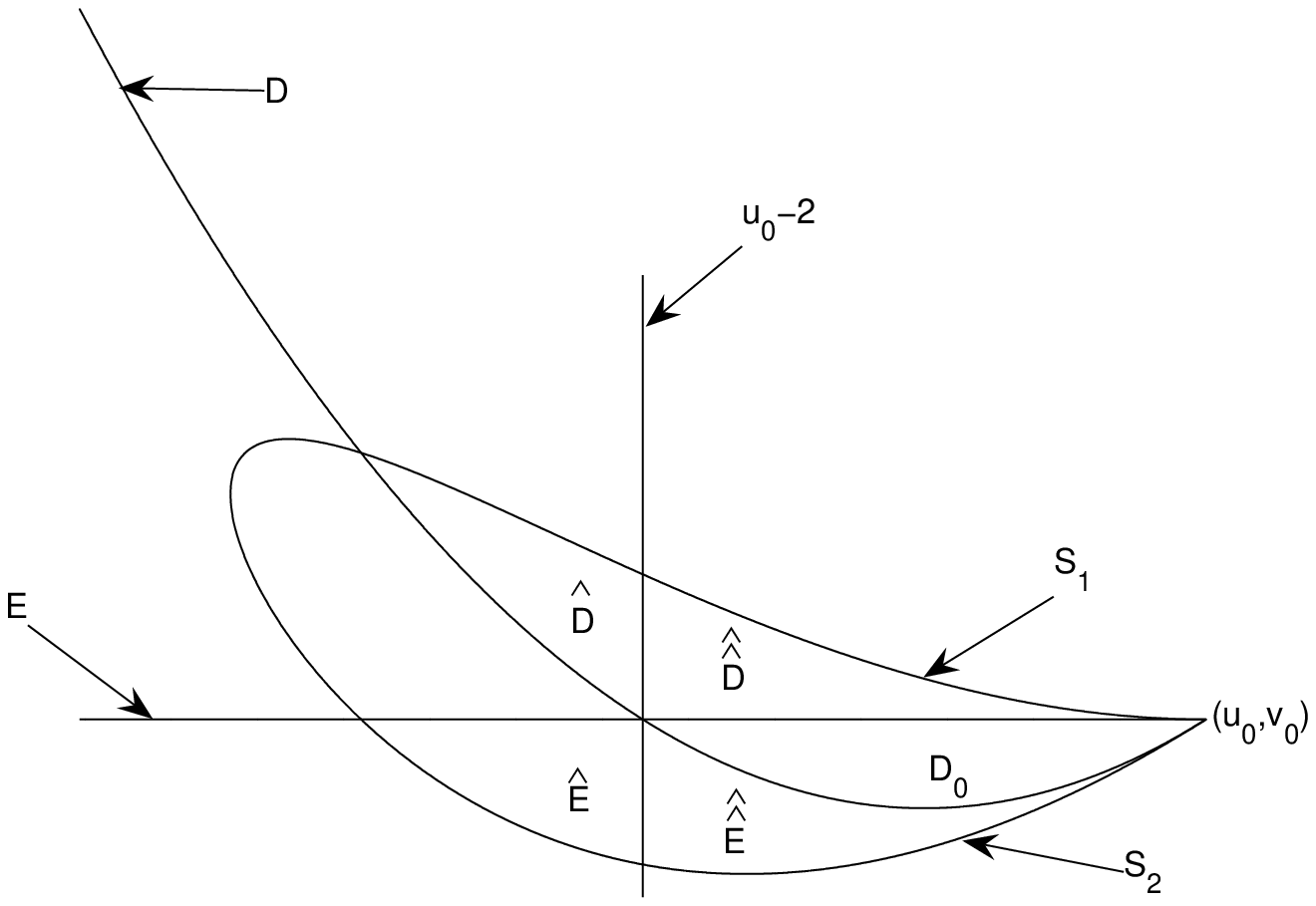}
%
\caption{}
\end{figure}

One can see that a rarefaction curve cannot enter into $D_{0}$
since it has to pass trough
the intersection point $(u_{0}-2,-2u_{0}+v_{0}+2)$ of $D$ and $E$, but 
\begin{equation*}
\begin{split}
& 2u-2u_{0}-uu_{0}+u^{2}/2+u_{0}^{2}/2+2>0
\text{ (ie. } R_{1} 
\text{ is above the curve  } E\text{)} \\
& 2u_{0}-2u+uu_{0}-u^{2}/2-u_{0}^{2}/2-2<0 \text{ (ie. } R_{2} 
\text{ is bellow the curve } D\text{)}. 
\end{split}
\end{equation*}
Therefore, a rarefaction curve which passes trough the point
$D\cap E$ goes either into $\hat{\hat{D}}$ or $\hat{\hat{E}}$,
and these cases are analysed above.

Thus, we have described all important points of the interactions
between singular shock and rarefaction waves. 
When a result of a single interaction is known, the question
about further singular shock path could be answered by a successive 
use of the above procedures.

\begin{remark} \label{r4} 
One can use the similar analysis of all possible cases when  
a rarefaction wave which interacts with a singular shock wave is on
the left-hand side of it. Instead of direct rarefaction and singular
shock curves, the inverse ones should be used, i.e.\
$(u_{2},v_{2})$ is a starting point and one is able to calculate 
$v_{0}$ from formulas of $E$ ,$D$, $S_{1}$, $S_{2}$, $R_{1}$ and $R_{2}$.
\end{remark} 

\begin{remark}
In the contrast with the case in \cite{NedObe}, where interaction 
can generate some ``strange'' solution containing unbounded $L_{loc}^{1}$
function, in the presented system one can find only bounded functions
and singular shock waves as a result on an interaction.

For a system (\ref{gdss1}) with $g_{1}\not \equiv \mathop{\rm const}$, or
$g_{2}\not \equiv 0$, interaction of singular shock and rarefaction waves
cannot be treated as easy as here. 
\end{remark}
   
Thus, we have proved the following assertion for the interaction
in the case of system (\ref{dss1}).

\begin{theorem} \label{primenjena}
Suppose that a singular shock wave interacts with a rarefaction wave
at the time $T$. For some time period $T<t<T_{1}$ the solution
is represented by a singular shock wave supported by a uniquely defined
curve (not a line) followed by a new rarefaction wave. Depending on the 
right-hand value of the primary rarefaction wave, one has the following 
possible cases for a solution after $t>T_{1}$.
\begin{itemize}
\item[(a)] Single singular shock wave (supported by a line)
with an increasing strength. 
\item[(b)] 1-rarefaction wave followed by singular shock wave
with an increasing strength. 
\item[(c)] Singular shock wave with an increasing strength
followed by 2-rarefaction wave.
\item[(d)] Singular shock wave with a decreasing strength 
prolonged by either a single singular shock wave with an increasing strength,
or a pair of admissible shock waves.
\item[(e)]  1-rarefaction wave followed by singular shock wave
with a constant strength. 
\item[(f)]  Singular shock wave with a constant strength
followed by 2-rarefaction wave.
\item[(g)] Singular shock wave with a decreasing strength
prolonged by either 1-rarefaction wave followed by singular shock wave
with decreasing or constant strength, or singular shock wave
with decreasing or constant strength followed by 2-rarefaction wave.
\end{itemize}
``Prolonged'' is the state after strength of singular shock wave 
becomes zero. Such wave can also stop with non-zero strength, and then
there is obviously no prolongation described above. 
\end{theorem}

\end{document}